\documentclass[11pt]{amsart}
\usepackage{graphicx} 

\usepackage{amsfonts}
\usepackage{amsmath}
\usepackage{amsthm}
\usepackage{latexsym}
\usepackage{amssymb}

\usepackage{amscd}
\usepackage{xcolor}
\usepackage{hyperref}
\usepackage{alphabeta}
\usepackage{fancyhdr}
\usepackage{mathtools}
\usepackage{fancyvrb}

\usepackage{imakeidx}
\usepackage{calrsfs}
\usepackage{soul}
\usepackage{multirow}
\usepackage{todonotes}
\usepackage{amsbsy}
\usepackage[all]{xy}
\usepackage{changes}
\usepackage{enumerate}
\usepackage{comment}
\usepackage{mathtools}
\usepackage{float}
\usepackage{color, colortbl}
\usepackage{changes}
\usepackage{makeidx}

\newtheorem{theorem}{Theorem}

\theoremstyle{definition}

\theoremstyle{remark}

\numberwithin{equation}{section}

\def\cH{{\mathcal H}}

\begin{document}

\title[A variant of Hilbert's inequality and the norm of $\cH$ on $K^p$]{A variant of Hilbert's inequality and the norm of the Hilbert Matrix on $K^p$}

\author{V. Daskalogiannis}
\email{vdaskalo@math.auth.gr}
\thanks{This first author was supported by the Hellenic Foundation for Research and Innovation (H.F.R.I.) under the ``2nd Call for H.F.R.I. Research Projects to support Faculty Members \& Researchers'' (Project Number: 73342).}

\author{P. Galanopoulos}
\email{petrosgala@math.auth.gr}
\address{Department of Mathematics, Aristotle University of Thessaloniki, 54124, Greece }

\author{M. Papadimitrakis}
\email{mpapadim@uoc.gr}
\address{Department of Mathematics and Applied Mathematics, University of Crete, 70013, Heraklion, Greece}

\begin{abstract}
We prove the nontrivial variant
	\[
	\sum\limits_{m,n=1}^{\infty}\Big(\frac{n}{m}\Big)^{\frac{1}{q}-\frac{1}{p}}\frac{a_mb_n}{m+n-1}\leq\frac{\pi}{\sin\frac{\pi}{p}} \Big( \sum\limits_{m=1}^{\infty}a_m^p\Big)^{\frac 1p}\Big( \sum\limits_{n=1}^{\infty}b_n^q\Big)^{\frac 1q}
	\]
of the well known Hilbert's inequality. Then we use this to determine the exact value $\frac{\pi}{\sin\frac{\pi}p}$ of the norm of the Hilbert matrix as an operator acting on the Hardy-Littlewood space $K^p$. This space consists of all functions $f(z)=\sum\limits_{m=0}^{\infty}a_mz^m$ analytic in the unit disc with $\|f\|_{K^p}^p=\sum\limits_{m=0}^{\infty}(m+1)^{p-2}|a_m|^p<+\infty$.
\end{abstract}

\maketitle

\section{Preliminaries}

The Hilbert matrix is the infinite matrix, whose entries are 
\[
 \frac{1}{m+n-1},\quad n,m=1,2\ldots.
\]

The well known Hilbert\rq{}s inequality \cite[Th. 323]{Hardy1934} states that if $(a_m),(b_n)$ are sequences of non negative terms such that $(a_m)\in\ell^p$, $(b_n)\in\ell^q$, then
\begin{equation}\label{hardy}
\sum_{m,n=1}^{\infty}\frac{a_mb_n}{m+n-1}\leq
\frac{\pi}{\sin\frac{\pi}{p}}\Big(\sum_{m=1}^{\infty}{a_m}^p\Big)^{\frac{1}{p}}\Big(\sum_{n=1}^{\infty}{b_n}^q\Big)^{\frac{1}{q}},
\end{equation}
where $1<p<\infty$, $1<q<\infty$, $\frac{1}{p}+\frac{1}{q}=1$, and the constant $\frac{\pi}{\sin\frac{\pi}{p}}$ is the smallest possible for this inequality. This implies that the Hilbert matrix induces a bounded operator $\mathcal{H}$,
\[
\mathcal{H}:(a_m)\longmapsto\,\mathcal{H}(a_m)=\Big(\sum_{m=1}^{\infty}\frac{a_m}{m+n-1}\Big)
\]
on the spaces $\ell^p$, $1<p<\infty$, with norm
\[
 \|\mathcal{H}\|_{\ell^p \to \ell^p}=\frac{\pi}{\sin\frac{\pi}{p}}.
\]

The operator $\mathcal{H}$ can also be considered as an operator on spaces of analytic functions by its action on the sequence of Taylor coefficients of any such function. 

Let $\mathbb{D}=\{z\in \mathbb{C}:\,|z|<1\}$ be the unit disk and  $H(\mathbb{D})$ be the space of analytic functions on $\mathbb{D}$. 

The Hardy space $H^p$, $0<p<\infty$, consists of all $f\in H(\mathbb{D})$ for which
\[
 \|f\|_{H^p}=\sup_{0\leq r<1}M_p(r,f)<\infty,
\]
where $M_p^p(r,f)$ are the integral means
\[
M^p_p(r,f)=\frac{1}{2\pi}\int_{0}^{2\pi}|f(re^{i\theta})|^p\,d\theta.
\]
If $p\geq 1$, then $H^p$ is a Banach space under the norm $\|\cdot\|_{H^p}$. If $0<p<1$, then $H^p$ is a complete metric space.

For $f(z)=\sum\limits_{m=0}^{\infty}a_m z^m\,\in {H}^1$, Hardy\rq{}s inequality \cite[p.48]{Duren1970}
\[
 \sum_{m=0}^\infty\frac{|a_m|}{m+1}\leq\pi\|f\|_{H^1},
\]
 implies that the power series
\[
  \mathcal{H}(f)(z)=\sum_{n=0}^\infty \Big(\sum_{m=0}^\infty \frac{a_m}{m+n+1}\Big)z^n
\]
has bounded coefficients. Therefore $\mathcal{H}(f)$ is an analytic function of the unit disk for any $f\in H^1$ and hence for any $f\in H^p$, $p\geq 1$.

The Bergman space $A^p$, $0<p<\infty$, consists of all $f\in H(\mathbb{D})$ for which
\[
 \|f\|_{A^p}=\Big(\int_\mathbb{D} |f(z)|^p\, dA(z)\Big)^{1/p}<\infty,
\]
where $dA(z)$ is the normalized Lebesgue area measure on $\mathbb{D}$.  
If $p\geq 1$, then $A^p$ is a Banach space under the norm $\|\cdot\|_{A^p}$.
 
If $f(z)=\sum\limits_{m=0}^{\infty}a_m z^m\in A^p$ and $p>2$, then by \cite[Lemma 4.1]{Nowak2010} we have
\[
\sum_{m=0}^\infty\frac{|a_m|}{m+1}<\infty.
\]
Thus $\mathcal{H}(f)$ is an analytic function in $\mathbb{D}$ for each function $f\in A^p$, $p>2$.

E. Diamantopoulos and A. G. Siskakis initiated the study of the Hilbert matrix as an operator on Hardy and Bergman spaces in \cite{Diamantopoulos2004, Diamantopoulos2000} and showed that $\mathcal{H}(f)$ has the following integral representation
\[
\mathcal{H}(f)(z)=\int_{0}^{1}\frac{f(t)}{1-tz}\,dt,\quad z\in\mathbb{D}.
\]
Then, considering $\mathcal{H}$ as an average of weighted composition operators, they showed that it is a bounded operator on $H^p$, $p>1$, and on $A^p$, $p>2$, and they estimated its norm. Their study was further extended by M. Dostani{\'{c}}, M. Jevti{\'{c}} and D. Vukoti{\'{c}} in \cite{Dostanic2008} and by V. Bo{\v{z}}in and B. Karapetrovi{\'{c}} in \cite{Bozin2018} (see also \cite{Lindstrom2018}). Summarizing their results, we now know that
\[
 \|\mathcal{H}\|_{H^p \to H^p}=\|\mathcal{H}\|_{A^{2p}\to A^{2p}}=\frac{\pi}{\sin\frac{\pi}{p}},\quad 1<p<\infty.
\]

We define $K^p$, $\,0<p<\infty$, to be the space of all $f(z)=\sum\limits_{m=0}^{\infty}a_mz^m \in H(\mathbb{D})$ such that 
\[
\|f\|^p_{K^p}=\sum_{m=0}^{\infty}(m+1)^{p-2}|a_m|^p<\infty.
\]
If $p\geq 1$, then $K^p$ is a Banach space under the norm $\|\cdot\|_{K^p}$.

According to the classical Hardy-Littlewood inequalities, \cite[Th. 6.2 \& 6.3]{Duren1970}, if 
$f(z)=\sum\limits_{m=0}^{\infty}a_mz^m\in H^p$, $0<p\leq 2$, then
\[
\sum_{m=0}^{\infty}(m+1)^{p-2}|a_m|^p\leq c_p \|f\|_{H^p}^p
\]
and hence $f\in K^p$. Also, if $2\leq p<\infty$ and $f(z)=\sum\limits_{m=0}^{\infty}a_mz^m\in K^p$, then
\[
\|f\|_{H^p}^p\leq c_p\sum_{m=0}^{\infty}(m+1)^{p-2}|a_m|^p
\]
and hence $f\in H^p$. In both cases $c_p$ is a constant independent of $f$.

If $p\geq 1$, and in the special case where the sequence $(a_m)$ is real and decreasing to zero, then for $f(z)=\sum\limits_{m=0}^{\infty}a_mz^m$ we have that $f\in H^p$ if and only if $f\in K^p$ \cite[Th. A \& 1.1]{pav1} .

Now it is clear that the proper domain of definition of the operator $\mathcal{H}$ acting on analytic functions in the unit disc is the space $K^1$. Indeed, if $f(z)=\sum\limits_{m=0}^{\infty} a_mz^m\,\in K^1$ then 
\[
\sum_{m=0}^\infty\frac{|a_m|}{m+1}<\infty,
\]
and hence $\mathcal{H}(f)\in H(\mathbb{D})$.

Moreover, when $1<p<\infty$ and $f\in K^p$, we consider $q$ so that $\frac{1}{p}+\frac{1}{q}=1$ and we apply H{\"o}lder's inequality to find
\[
\begin{split}
\sum_{m=0}^\infty\frac{|a_m|}{m+1}&=\sum_{m=0}^\infty (m+1)^{\frac{2}{p}-2}(m+1)^{1-\frac{2}{p}}|a_m| \\
 &\leq\Big(\sum_{m=0}^\infty\frac{1}{(m+1)^{2}} \Big)^{\frac{1}{q}}\Big(\sum_{m=0}^\infty (m+1)^{p-2}|a_m|^p \Big)^{\frac{1}{p}}<\infty.
\end{split}
\]
Hence $K^p\subseteq K^1$ and so, if $f\in K^p$, then $\mathcal{H}(f)$ defines an analytic function in $\mathbb{D}$.

Recently, in \cite[Theorem 1]{Pelaez2022} (see also \cite{Contreras2016}), the authors associated the boundedness of the generalized Volterra operators 
\[
T_g(z)=\int_0^z f(w) g'(w) dw,\quad z\in\mathbb{D},
\] 
induced by symbols $g\in H(\mathbb{D})$ with non-negative Taylor coefficients and acting from a space $X$ to $H^\infty$, to the $K^p$-norm of the function $\mathcal{H}(g')$. In this result $X$ can be $H^p$ or $K^p$ or the Dirichlet-type space $D^p_{p-1}$.

\section{A variant of Hilbert's inequality}

Our first result is a nontrivial variant of the classical Hilbert's inequality.

Before we state our first main result we shall mention two more variants of Hilbert's inequality. The first, in \cite{Yang2004}, is
\[
\sum\limits_{m,n=1}^{\infty}\Big(\frac{n}{m}\Big)^{\frac{1}{q}-\frac{1}{p}}\frac{a_mb_n}{m+n}\leq\frac{\pi}{\sin\frac{\pi}{p}} \Big( \sum\limits_{m=1}^{\infty}a_m^p\Big)^{\frac 1p}\Big( \sum\limits_{n=1}^{\infty}b_n^q\Big)^{\frac 1q}
\]
and the second, in \cite{Yang2005}, is
\[
\sum\limits_{m,n=1}^{\infty}\Big(\frac{n-\frac 12}{m-\frac 12}\Big)^{\frac{1}{q}-\frac{1}{p}}\frac{a_mb_n}{m+n-1}\leq\frac{\pi}{\sin\frac{\pi}{p}} \Big( \sum\limits_{m=1}^{\infty}a_m^p\Big)^{\frac 1p}\Big( \sum\limits_{n=1}^{\infty}b_n^q\Big)^{\frac 1q}.
\]
In fact Yang proves a whole family of such inequalities depending on a parameter. In all these variants, as well as in the original Hilbert's inequality, the kernel involved in the double sum is of the form
\[
\Big(\frac{k(n)}{k(m)}\Big)^{c_p}\frac 1{k(m)+k(n)}
\]
which is homogeneous of degree $-1$. As a consequence, in order to prove these variants one needs to apply the standard arguments used in the proof of the original Hilbert's inequality. The kernel in our variant of Hilbert's inequality lacks this homegeneity and the standard arguments do not apply.
 
\begin{theorem}\label{theorem}\label{new ineq}
Let $1<p<\infty$, $1<q<\infty$, $\frac{1}{p}+\frac{1}{q}=1$. If $(a_m)\in \ell^p$, $(b_n)\in \ell^q$ are sequences of non negative terms, then
\[
\sum\limits_{m,n=1}^{\infty}\Big(\frac{n}{m}\Big)^{\frac{1}{q}-\frac{1}{p}}\frac{a_mb_n}{m+n-1}\leq\frac{\pi}{\sin\frac{\pi}{p}} \Big( \sum\limits_{m=1}^{\infty}a_m^p\Big)^{\frac 1p}\Big( \sum\limits_{n=1}^{\infty}b_n^q\Big)^{\frac 1q}.
\]
The constant $\frac{\pi}{\sin\frac{\pi}{p}}$ is the smallest possible for this inequality.
\end{theorem}
\begin{proof}
In fact we may restrict to $1<q\leq 2\leq p<\infty$.\newline
We assume
\[
\frac{\alpha}{p}+\frac{\beta}{q}=1,\quad\alpha\geq 0,\,\beta\geq 0,
\]
where $\alpha$ and $\beta$ will be chosen appropriately later.\newline 
By H{\"o}lder's inequality,
\begin{align*}
	&\sum\limits_{m,n=1}^{\infty}\Big(\frac{n}{m}\Big)^{\frac{1}{q}-\frac{1}{p}}\frac{a_mb_n}{m+n-1} \\
	&=\sum\limits_{m,n=1}^{\infty}\Big(\frac{n}{m}\Big)^{(\frac{1}{pq}-\frac{1}{p})+(\frac{1}{q}-\frac{1}{pq})}\frac{a_mb_n}{(m+n)^{\frac{1}{p}}(m+n)^{\frac{1}{q}}}\Big(\frac{m+n}{m+n-1}\Big)^{\frac{\alpha}{p}}\Big(\frac{m+n}{m+n-1}\Big)^{\frac{\beta}{q}} \\
	&\leq \Big(\sum\limits_{m=1}^{\infty} a_m^p \Big(\sum\limits_{n=1}^{\infty}\Big(\frac{m}{n} \Big)^{\frac{1}{p}}\frac{1}{(m+n)^{1-\alpha}(m+n-1)^{\alpha}}\Big)\Big)^{\frac{1}{p}}  \\
	&\times\Big(\sum\limits_{n=1}^{\infty} b_n^q\Big(\sum\limits_{m=1}^{\infty}\Big(\frac{n}{m} \Big)^{\frac{1}{q}}\frac{1}{(m+n)^{1-\beta}(m+n-1)^{\beta}}\Big)\Big)^{\frac{1}{q}}.
\end{align*}
Hence it is enough to prove
\begin{align}\label{ineq 1}
	\sum\limits_{n=1}^{\infty}\Big(\frac{m}{n} \Big)^{\frac{1}{p}}\frac{1}{(m+n)^{1-\alpha}(m+n-1)^{\alpha}}\leq\frac{\pi}{\sin\frac{\pi}{p}},\quad m\geq 1,
\end{align}
and
\begin{align}\label{ineq 2}
\sum\limits_{m=1}^{\infty}\Big(\frac{n}{m} \Big)^{\frac{1}{q}}\frac{1}{(m+n)^{1-\beta}(m+n-1)^{\beta}}\leq\frac{\pi}{\sin\frac{\pi}{q}},\quad n\geq 1,
\end{align}
where, of course, $\sin\frac{\pi}{p}=\sin\frac{\pi}{q}$.\newline
Now we observe that, for all $\alpha\geq 0$, $p>0$, $m\geq 1$, the positive function
\[ 
 f(t)=t^{-\frac{1}{p}}(m+t)^{\alpha-1}(m+t-1)^{-\alpha}, \quad t>0,
\]
is convex. Indeed, taking the second derivative of the logarithm of $f(t)$, we get
\[
 \frac{f(t)f''(t)-f'(t)^2}{f(t)^2}=\frac{t^{-2}}{p}+ (m+t)^{-2}+\alpha\big((m+t-1)^{-2}-(m+t)^{-2}\big)>0,
\]
which proves that $f''(t)>0$. In fact, this calculation proves more: that $f$ is logarithmically convex.\newline
The convexity of $f$ implies
\[
 f(n)\leq\int_{n-\frac{1}{2}}^{n+\frac{1}{2}}f(t)\,dt,\quad n\geq 1.
\]
Adding these inequalities we get for the left side of \eqref{ineq 1} that
\begin{align*}
	\sum\limits_{n=1}^{\infty}\Big(\frac{m}{n} \Big)^{\frac{1}{p}}&\frac{1}{(m+n)^{1-\alpha}(m+n-1)^{\alpha}}\\
	&\leq \int_{\frac{1}{2}}^{\infty}\Big(\frac{m}{t} \Big)^{\frac{1}{p}}\frac{1}{(m+t)^{1-\alpha}(m+t-1)^{\alpha}}\,dt\\
	&=\int_{\frac{1}{2m}}^{\infty}\frac{1}{t^{\frac{1}{p}}(t+1)^{1-\alpha}(t+1-\frac{1}{m})^{\alpha}}\,dt
\end{align*}
by the change of variables $t\mapsto mt$.\newline
Therefore, in order to prove \eqref{ineq 1} it is enough to prove
\begin{equation}\label{ineq 3}
	\int_{\frac{1}{2m}}^{\infty}\frac{1}{t^{\frac{1}{p}}(t+1)^{1-\alpha}(t+1-\frac{1}{m})^{\alpha}}\,dt\leq\frac{\pi}{\sin\frac{\pi}{p}},\quad m\geq 1.
\end{equation}
We consider now the function
\begin{align*}
	F(y)&=\int_{y}^{\infty}\frac{1}{t^{\frac{1}{p}}(t+1)^{1-\alpha}(t+1-2y)^{\alpha}}\,dt \\
	&=\int_{0}^{\infty}\frac{1}{(t+y)^{\frac{1}{p}}(t+1+y)^{1-\alpha}(t+1-y)^{\alpha}}\,dt,\quad 0\leq y\leq\frac{1}{2}.
\end{align*}
Hence in order to prove \eqref{ineq 3} it is enough to prove 
\begin{equation}\label{ineq 4}
	F(y)\leq\frac{\pi}{\sin\frac{\pi}{p}},\quad 0\leq y\leq\frac{1}{2}.
\end{equation}
Now, exactly as before, we observe that, for all $\alpha\geq 0$, $p>0$, $t>0$, the positive function 
\[
 g_t(y)=(t+y)^{-\frac{1}{p}}(t+1+y)^{\alpha-1}(t+1-y)^{-\alpha},\quad 0\leq y\leq\frac{1}{2},
\]
is convex. Indeed, we take the second derivative of the logarithm of $g_t(y)$ and we get
\begin{align*}
	\frac{g_t(y)g_t''(y)-g_t'(y)^2}{g_t(y)^2}=&\frac{(t+y)^{-2}}{p}+(t+1+y)^{-2}\\
	&+\alpha\big((t+1-y)^{-2}-(t+1+y)^{-2}\big)>0,
\end{align*}
which proves that $g_t''(y)>0$.\newline
Thus $F(y)=\int_{0}^{\infty}g_t(y)\,dt$ is also convex and, as such, it satisfies
\[ 
 F(y)\leq \max\Big\{F(0), F\Big(\frac{1}{2}\Big)\Big\}.
\]
Since 
\[
 F(0)=\int_{0}^{\infty}\frac{1}{t^{\frac{1}{p}}(t+1)}\,dt=\frac{\pi}{\sin \frac{\pi}{p}},
\]
in order to prove \eqref{ineq 4} it is enough to prove 
\[
 F\Big(\frac{1}{2}\Big)\leq\frac{\pi}{\sin \frac{\pi}{p}}.
\]
Since
\[
 F\Big(\frac{1}{2}\Big)=\int_{1/2}^{\infty}\frac{\left(t+1\right)^{\alpha-1}}{t^{\frac{1}{p}+\alpha}}\,dt
	=\int_{0}^{2} \frac{(t+1)^{\alpha}}{t^{1-\frac{1}{p}}(t+1)}\,dt
\]
after the change of variables $t\mapsto 1/t$, we conclude that in order to prove \eqref{ineq 1} it is enough to prove
\[
 \int_{0}^{2} \frac{(t+1)^{\alpha}}{t^{1-\frac{1}{p}}(t+1)}\,dt\leq \frac{\pi}{\sin \frac{\pi}{p}}.
\]
In exactly the same manner, we see that in order to prove \eqref{ineq 2} it is enough to prove
\[
 \int_{0}^{2} \frac{(t+1)^{\beta}}{t^{1-\frac{1}{q}}(t+1)}\,dt\leq \frac{\pi}{\sin \frac{\pi}{q}}.
\]
We make the change of notation
\[
 x=\frac{1}{p}, \quad 1-x=\frac{1}{q},
\]
and, after $\frac{\alpha}{p}+\frac{\beta}{q}=1$, we write
\[
 \beta=\frac{1-\alpha x}{1-x},
\]
where $0\leq \alpha x\leq 1$. Then our last two inequalities become
\begin{equation}\label{ineq 5}	
	\int_{0}^{2}\frac{(t+1)^{\alpha}}{t^{1-x}(t+1)}\,dt\leq \frac{\pi}{\sin\pi x}=\int_0^{\infty}\frac{1}{t^{1-x}(t+1)}\,dt 
\end{equation} 
and
\begin{equation}\label{ineq 6}	
	\int_{0}^{2} \frac{(t+1)^{\frac{1-\alpha x}{1-x}}} {t^{x}(t+1)}\,dt\leq \frac{\pi}{\sin\pi x}=\int_0^{\infty}\frac{1}{t^{x}(t+1)}\,dt.
\end{equation}
Now, inequality \eqref{ineq 5} is equivalent to 
\[
 \int_{0}^{2} \frac{(t+1)^{\alpha}-1}{t^{1-x}(t+1)}\,dt\leq \int_{2}^{\infty}\frac{1}{t^{1-x}(t+1)}\,dt  
\]
or, after the change of variables $t\mapsto 2t$, to 
\[
 \int_{0}^{1}\frac{(2t+1)^{\alpha}-1}{t^{1-x}(2t+1)}\,dt\leq\int_{1}^{\infty}\frac{1}{t^{1-x}(2t+1)}\,dt,
\] 
or finally, substituting $t\mapsto 1/t$ in the left integral, to the inequality
\begin{equation}\label{ineqI}
	\int_{1}^{\infty} 	\frac{\left(1+\frac{2}{t}\right)^{\alpha}-1}{t^{x}(t+2)}\,dt\leq\int_{1}^{\infty}\frac{1}{t^{1-x}(2t+1)}\,dt, \quad 0<x\leq \frac{1}{2}.
\end{equation}
Similarly, inequality \eqref{ineq 6} is equivalent to
\[
 \int_{0}^{2} \frac{(t+1)^{\frac{1-\alpha x}{1-x}}-1}{t^{x}(t+1)}\,dt\leq \int_{2}^{\infty}\frac{1}{t^{x}(t+1)}\,dt
\]
or, after the successive change of variables $t\mapsto 2t$ and $t\mapsto 1/t$, to 
\begin{equation}\label{ineqII}
	\int_{1}^{\infty}\frac{\left(1+\frac{2}{t}\right)^{\frac{1-\alpha x}{1-x}}-1}{t^{1-x}(t+2)}\,dt\leq\int_{1}^{\infty}\frac{1}{t^{x}(2t+1)}\,dt, \quad 0<x\leq \frac{1}{2}.
\end{equation}
So we have come to the point where, for every $x$ with $0<x\leq\frac 12$, we have to prove inequalities \eqref{ineqI} and \eqref{ineqII} for a proper choice of $\alpha$ with $0\leq\alpha\leq\frac 1x$.\newline
A very usefull observation for what follows is that for fixed $\alpha$ with $0\leq\alpha\leq 1$, if \eqref{ineqI} holds for some $x$, then it holds for all larger $x$. The reason is that the left-hand side in \eqref{ineqI} is a decreasing function of $x$ and the right-hand side in \eqref{ineqI} is an increasing function of $x$. Similarly, if \eqref{ineqII} holds for some $x$, then it holds for all smaller $x$. It helps to see that for fixed $\alpha$ with $0\leq\alpha\leq 1$ the function $\frac{1-\alpha x}{1-x}$ is increasing.\newline
Now we split the interval $0<x\leq\frac 12$ in three subintervals in each of which we make the corresponding choices $\alpha=0$, $\alpha=1$ and $\alpha=\frac 12$.\vspace{0.2em}\newline
\textit{The case $\alpha=0$}.\vspace{0.2em}\newline 
Let $\alpha=0$. First of all, it is obvious that \eqref{ineqI} is true for all $0<x\leq\frac 12$. We claim that \eqref{ineqII} is valid for all $0<x\leq\frac 13$ and as we observed it is enough to prove it for $x=\frac 13$.\newline
Observe now that $0<x\leq\frac{1}{2}$ implies 
$0<\frac{x}{1-x}\leq 1,$
so by Bernoulli's inequality we get
\begin{align*}
	\Big(1+\frac{2}{t}\Big)^{\frac{1}{1-x}}&=\Big(1+\frac{2}{t}\Big)\Big(1+\frac{2}{t}\Big)^{\frac{x}{1-x}}\leq\Big(1+\frac{2}{t}\Big)\Big(1+\frac{x}{1-x}\,\frac{2}{t}\Big)\\
	&=1+\frac{2}{t}+\frac{x}{1-x}\,\frac{2(t+2)}{t^2}.
\end{align*}
Hence 
\[
 \int_{1}^{\infty} \frac{\left(1+\frac{2}{t}\right)^{\frac{1}{1-x}}-1} {t^{1-x}(t+2)}\,dt \leq \int_{1}^{\infty}\frac{2}{t^{2-x}(t+2)}\,dt+\frac{2x}{1-x}\int_{1}^{\infty}\frac{1}{t^{3-x}}\,dt.
\]
Using 
\begin{equation}\label{split}
	\frac{2}{t(t+2)}=\frac{1}{t}-\frac{1}{t+2}
\end{equation} 
the last inequality becomes
\begin{align*}
	\int_{1}^{\infty} \frac{\left(1+\frac{2}{t}\right)^{\frac{1}{1-x}}-1} {t^{1-x}(t+2)}\,dt &\leq
	\int_{1}^{\infty}\frac{1}{t^{2-x}}\,dt-\int_{1}^{\infty}\frac{1}{t^{1-x}(t+2)}\,dt+\frac{2x}{(1-x)(2-x)}\\
	&=\frac{2+x}{(1-x)(2-x)}-\int_{1}^{\infty}\frac{1}{t^{1-x}(t+2)}\,dt.
\end{align*}
Hence in order to prove \eqref{ineqII} we need to have  
\begin{align*}
	\frac{2+x}{(1-x)(2-x)}&\leq \int_{1}^{\infty} \frac{1}{t^{1-x}(t+2)}\,dt+\int_{1}^{\infty} \frac{1}{t^{x}(2t+1)}\,dt\\
	&=\int_{0}^{1} \frac{1}{t^{x}(2t+1)}\,dt+\int_{1}^{\infty} \frac{1}{t^{x}(2t+1)}\,dt=\int_{0}^{\infty}\frac{1}{t^{x}(2t+1)}\,dt\\
	&=2^{x-1}\int_{0}^{\infty}\frac{1}{t^{x}(t+1)}\,dt=2^{x-1}\frac{\pi}{\sin\pi x}.
\end{align*}
For $x=\frac 13$ this becomes $\frac{21}{10}\leq \frac{2^{\frac 13}\pi}{\sqrt{3}}$ which is true and proves our claim.\newline
We proved that when $\alpha=0$ both \eqref{ineqI} and \eqref{ineqII} hold for $0<x\leq\frac 13$.\vspace{0.2em}\newline
\textit{The case $\alpha=1$}.\vspace{0.2em}\newline 
Let $\alpha=1$. In this case \eqref{ineqI} becomes 
\begin{equation}\label{ineqI a=1}
	\int_{1}^{\infty}\frac{2}{t^{1+x}(t+2)}\,dt\leq\int_{1}^{\infty}\frac{1}{t^{1-x}(2t+1)}\,dt.
\end{equation}
We claim that this inequality is true for $\frac 25\leq x\leq\frac 12$ and it suffices to prove it for $x=\frac 25$.\newline
Using \eqref{split}, the left-hand side of \eqref{ineqI a=1} becomes
\begin{align*}
	\int_{1}^{\infty}\frac{2}{t^{1+x}(t+2)}\,dt&=\int_{1}^{\infty} \frac{1}{t^{1+x}}\,dt-\int_{1}^{\infty} \frac{1}{t^{x}(t+2)}\,dt\notag \\
	&=\frac{1}{x}-\int_{1}^{\infty} \frac{1}{t^{x}(t+2)}\,dt,
\end{align*}
Therefore, \eqref{ineqI a=1} amounts to showing 
\begin{align*}
	\frac{1}{x}&\leq \int_{1}^{\infty}\frac{1}{t^x(t+2)}\,dt+\int_{1}^{\infty}\frac{1}{t^{1-x}(2t+1)}\,dt=\int_{0}^{\infty}\frac{1}{t^x(t+2)}\,dt\\
	&=2^{-x}\int_{0}^{\infty}\frac{1}{t^x(t+1)}\,dt=2^{-x}\,\frac{\pi}{\sin (\pi x)}.
\end{align*}
for $x=\frac 25$. Equivalently, we need to show that
$$\frac{\sin \pi x}{\pi x}\leq 2^{-x}$$
for $x=\frac 25$. Indeed we have that
$$\frac{\sin\frac{2\pi}5}{\frac{2\pi}5}<1-\frac 1{3!}\Big(\frac{2\pi}5\Big)^2+\frac 1{5!}\Big(\frac{2\pi}5\Big)^4<2^{-\frac 25}$$
as we easily see after a few calculations.\newline 
Thus, \eqref{ineqI} is valid for $\frac 25\leq x\leq\frac 12$.\newline
We now turn to \eqref{ineqII}, and we claim that it holds for $0<x\leq\frac 12$ and it suffices to prove it for $x=\frac 12$.\newline 
When $\alpha=1$, \eqref{ineqII} becomes
\[ 
 \int_1^{\infty}\frac{2}{t^{2-x}(t+2)}\,dt\leq \int_1^{\infty}\frac{1}{t^{x}(2t+1)}\,dt
\]
or, by the use of \eqref{split},
\[
 \int_1^{\infty}\frac{1}{t^{2-x}}\,dt- \int_1^{\infty}\frac{1}{t^{1-x}(t+2)}\,dt\leq \int_1^{\infty}\frac{1}{t^{x}(2t+1)}\,dt.
\]
This is equivalent to
\begin{align*}
\frac{1}{1-x}&\leq\int_1^{\infty}\frac{1}{t^{1-x}(t+2)}\,dt+\int_1^{\infty}\frac{1}{t^{x}(2t+1)}\,dt= \int_{0}^{\infty}\frac{1}{t^{1-x}(t+2)}\,dt\\
&=2^{x-1}\frac{\pi}{\sin \pi x}.
\end{align*}
When $x=\frac 12$ this becomes $2\sqrt{2}\leq \pi$ and it is clearly true.\newline
We proved that when $\alpha=1$ both \eqref{ineqI} and \eqref{ineqII} hold for $\frac 25\leq x\leq\frac 12$.\vspace{0.2em}\newline
\textit{The case $\alpha=\frac 12$}.\vspace{0.2em}\newline 
Let $\alpha=\frac 12$. We first deal with inequality \eqref{ineqI}, which we shall prove for $\frac 13\leq x\leq\frac 25$. As we know it is enough to prove it for $x=\frac 13$.\newline 
When $\alpha=\frac 12$, \eqref{ineqI} becomes 
\[
 \int_{1}^{\infty} \frac{\left(1+\frac{2}{t}\right)^{\frac 12}-1} {t^{x}(t+2)}\,dt \leq \int_{1}^{\infty}\frac{1}{t^{1-x}(2t+1)}\,dt.
\]
Bernoulli's inequality gives
\[
 \Big(1+\frac 2t\Big)^{\frac 12}\leq 1+\frac 12\,\frac 2t=1+\frac 1t
\]
and hence 
\[
 \int_{1}^{\infty} \frac{\left(1+\frac{2}{t}\right)^{\frac 12}-1} {t^{x}(t+2)}\,dt \leq \int_{1}^{\infty} \frac{1} {t^{1+x}(t+2)}\,dt.
\]
Therefore it suffices to show that 
\[
 \int_{1}^{\infty} \frac{1} {t^{1+x}(t+2)}\,dt \leq \int_{1}^{\infty}\frac{1}{t^{1-x}(2t+1)}\,dt
\]
for $x=\frac 13$. This is indeed true, since
\[ 
 t^{\frac 23}(2t+1)\leq t^{\frac 43}(t+2),\quad t\geq 1,
\]
as we easily see by raising to the third power.\newline
We now turn to \eqref{ineqII} which for $\alpha=\frac 12$ becomes 
\[
 \int_{1}^{\infty} \frac{\left(1+\frac{2}{t}\right)^{\frac{1}{2}\frac{2-x}{1-x}}-1} {t^{1-x}(t+2)}\,dt \leq \int_{1}^{\infty}\frac{1}{t^{x}(2t+1)}\,dt,
\]
and we claim it holds for $\frac 13\leq x\leq\frac 25$. Again it suffices to prove this inequality for $x=\frac 25$. Namely, it suffices to show 
\begin{equation}\label{ineqII a=1/2}
	\int_{1}^{\infty} \frac{\left(1+\frac{2}{t}\right)^{\frac 43}-1} {t^{\frac 35}(t+2)}\,dt \leq \int_{1}^{\infty}\frac{1}{t^{\frac 25}(2t+1)}\,dt.
\end{equation}
Taking into account Bernoulli's inequality, we have
\[
\Big(1+\frac{2}{t}\Big)^{\frac 43}= \Big(1+\frac{2}{t}\Big)\Big(1+\frac{2}{t}\Big)^{\frac 13} \leq\Big(1+\frac{2}{t}\Big)\Big(1+\frac 13\,\frac{2}{t}\Big)=1+\frac{4}{3t^2}(2t+1), 
\]
so instead of \eqref{ineqII a=1/2}, it suffices to prove
\begin{align}\label{ineqII a=1/2 new}
	\frac{4}{3}\int_{1}^{\infty} \frac{2t+1}{t^{2+\frac{3}{5}}(t+2)}\,dt\leq \int_{1}^{\infty}\frac{1}{t^{\frac{2}{5}}(2t+1)}\,dt.
\end{align}
Observe that the left-hand side of \eqref{ineqII a=1/2 new}, in view of \eqref{split}, is equal to 
\begin{align*}
	\frac{4}{3}& \int_{1}^{\infty}\frac{2t+1}{t^{2+\frac{3}{5}}(t+2)}\,dt = \frac{2}{3} \int_{1}^{\infty} \frac{2t+1}{t^{2+\frac{3}{5}}}\,dt - \frac{2}{3} \int_{1}^{\infty} \frac{2t+1}{t^{1+\frac{3}{5}}(t+2)}\,dt\\
	&=\frac{4}{3} \int_{1}^{\infty} \frac{1}{t^{1+\frac{3}{5}}}\,dt+
	\frac{2}{3} \int_{1}^{\infty} \frac{1}{t^{2+\frac{3}{5}}}\,dt -\frac{4}{3} \int_{1}^{\infty} \frac{1}{t^{\frac{3}{5}}(t+2)}\,dt\\
	&-\frac{2}{3} \int_{1}^{\infty} \frac{1}{t^{1+\frac{3}{5}}(t+2)}\,dt\\
	&=\frac{20}{9}+\frac{5}{12}-\frac{4}{3} \int_{1}^{\infty} \frac{1}{t^{\frac{3}{5}}(t+2)}\,dt-\frac{1}{3} \int_{1}^{\infty} \frac{1}{t^{1+\frac{3}{5}}}\,dt+\frac{1}{3} \int_{1}^{\infty} \frac{1}{t^{\frac{3}{5}}(t+2)}\,dt,
\end{align*}
where we used \eqref{split} for the last equality. Thus, altogether we have
\[
 \frac{4}{3} \int_{1}^{\infty} \frac{2t+1}{t^{2+\frac{3}{5}}(t+2)}\,dt =\frac{25}{12}-\int_{1}^{\infty} \frac{1}{t^{\frac{3}{5}}(t+2)}\,dt.
\]
Therefore, \eqref{ineqII a=1/2 new} is equivalent to the inequality
\[
 \frac{25}{12}\leq \int_{1}^{\infty} \frac{1}{t^{\frac{3}{5}}(t+2)}\,dt+\int_{1}^{\infty} \frac{1}{t^{\frac{2}{5}}(2t+1)}\,dt=\int_{0}^{\infty} \frac{1}{t^{\frac{2}{5}}(2t+1)}\,dt=\frac{2^{-\frac{3}{5}}\pi}{\sin\frac{3\pi}{5} }
\]
This inequality is an easy consequence of the inequality $\frac{\sin\frac{2\pi}5}{\frac{2\pi}5}<2^{-\frac 25}$ which we proved when we considered the case $\alpha=1$. Indeed
\[
\sin\frac{3\pi}{5}=\sin\frac{2\pi}5<\frac{2\pi}52^{-\frac 25}=\frac{2\pi}52^{-\frac 35}2^{\frac 15}<\frac{2\pi}52^{-\frac 35}\Big(1+\frac 15\Big)=\frac{12\pi}{25}2^{-\frac 35}.
\]
\newline
We proved that when $\alpha=\frac 12$ both \eqref{ineqI} and \eqref{ineqII} hold for $\frac 13\leq x\leq\frac 25$.\newline
We have proved the inequality of our theorem and now we shall show that the constant $\frac{\pi}{\sin\frac{\pi}{p}}$ is the best possible in this inequality. The proof follows the lines of Hardy's corresponding proof for the original Hilbert's inequality \cite[proof of Theorem 317, p. 232]{Hardy1934}, adapted to our weighted setting. For the sake of completeness, we provide the details.\newline
We consider any $\epsilon>0$ and the sequences $(a_m(\epsilon))$ and $(b_n(\epsilon))$ defined by
\[
a_m(\epsilon)=m^{-\frac{1+\epsilon}{p}},\quad b_n(\epsilon)=n^{-\frac{1+\epsilon}{q}}.
\]
We then have
\[
\|(a_m(\epsilon))\|^p_{\ell^p}=\sum_{m=1}^\infty\tfrac{1}{m^{1+\epsilon}}.
\]
Now, since $\frac{1}{x^{1+\epsilon}}$ is decreasing for $x\geq 1$, we have
\[
\frac{1}{\epsilon}=\int_1^\infty\frac{1}{x^{1+\epsilon}}\,dx\leq\sum_{m=1}^\infty\frac{1}{m^{1+\epsilon}}\leq 1+\int_1^\infty\frac{1}{x^{1+\epsilon}}\,dx=1+\frac{1}{\epsilon}.
\]
Setting $\phi(\epsilon)=\sum\limits_{m=1}^\infty\frac{1}{m^{1+\epsilon}}-\frac{1}{\epsilon}$, we get
\begin{equation}\label{eq 13}
\|(a_m(\epsilon))\|^p_{\ell^p}=\frac{1}{\epsilon}+\phi(\epsilon),\quad 0\leq \phi(\epsilon)\leq 1.
\end{equation}
 Respectively, setting $\psi(\epsilon)=\sum\limits_{n=1}^\infty\frac{1}{n^{1+\epsilon}}-\frac{1}{\epsilon}$, we have
\begin{equation}\label{eq 14}
\|(b_n(\epsilon)\|^q_{\ell^q}=\frac{1}{\epsilon}+\psi(\epsilon),\quad 0\leq \psi(\epsilon)\leq 1.
\end{equation}
In addition, we have that
\begin{equation}\label{ineq 15}
\sum_{m,n=1}^\infty\Big(\frac{n}{m}\Big)^{\frac{1}{q}-\frac{1}{p}}\frac{a_m(\epsilon)b_n(\epsilon)}{m+n-1}\geq
\sum_{m,n=1}^\infty\Big(\frac{n}{m}\Big)^{\frac{1}{q}-\frac{1}{p}}\frac{a_m(\epsilon)b_n(\epsilon)}{m+n}.
\end{equation}
Now for $(x,y)$ in the square $[m,m+1)\times[n,n+1)$, $m\geq 1$, $n\geq 1$, we have
\begin{align*}
	\Big(\frac{n}{m}\Big)^{\frac{1}{q}-\frac{1}{p}}\frac{a_m(\epsilon)b_n(\epsilon)}{m+n}&=\Big(\frac{n}{m}\Big)^{\frac{1}{q}-\frac{1}{p}}
	\frac{m^{-\frac{1+\epsilon}{p}}n^{-\frac{1+\epsilon}{q}}}{m+n}=\frac{m^{-\frac{1}{q}-\frac{\epsilon}{p}}n^{-\frac{1}{p}-\frac{\epsilon}{q}}}{m+n}\\
	&\geq \frac{x^{-\frac{1}{q}-\frac{\epsilon}{p}}y^{-\frac{1}{p}-\frac{\epsilon}{q}}}{x+y}=\Big(\frac{y}{x}\Big)^{\frac{1}{q}-\frac{1}{p}}
	\frac{x^{-\frac{1+\epsilon}{p}}y^{-\frac{1+\epsilon}{q}}}{x+y}.
\end{align*}
Therefore
\begin{equation}\label{ineq 16}
\sum_{m,n=1}^\infty\Big(\frac{n}{m}\Big)^{\frac{1}{q}-\frac{1}{p}}\frac{a_m(\epsilon)b_n(\epsilon)}{m+n}\geq I(\epsilon),
\end{equation}
where $I(\epsilon)$ is defined by
\[
I(\epsilon)=\int_1^\infty \int_1^\infty
 \Big(\frac{y}{x}\Big)^{\frac{1}{q}-\frac{1}{p}}
\frac{x^{-\frac{1+\epsilon}{p}}y^{-\frac{1+\epsilon}{q}}}{x+y}\,dx\,dy
=\int_1^\infty \int_1^\infty
\frac{x^{-\frac{1}{q}-\frac{\epsilon}{p}}y^{-\frac{1}{p}-\frac{\epsilon}{q}}}{x+y}\,dx\,dy.
\]
Applying the change of variables $y\mapsto xy$, we get
\[
I(\epsilon)=\int_1^\infty\frac{1}{x^{1+\epsilon}}\int_{\frac{1}{x}}^\infty \frac{1}{y^{\frac{1}{p}+\frac{\epsilon}{q}}(1+y)}\,dy\,dx
\]
Another change of variables $x\mapsto \frac{1}{x}$ gives
\begin{align*}
I(\epsilon)&=\int_0^1 x^{\epsilon-1}\int_x^\infty \frac{1}{y^{\frac{1}{p}+\frac{\epsilon}{q}}(1+y)}\,dy\,dx
\\
&=\int_0^1\frac{1}{\epsilon}(x^\epsilon)'\int_x^\infty\frac{1}{y^{\frac{1}{p}+\frac{\epsilon}{q}}(1+y)}\,dy\,dx\\
&=\frac{1}{\epsilon}\Big(\int_1^\infty \frac{1}{y^{\frac{1}{p}+\frac{\epsilon}{q}}(1+y)}\,dy+\int_0^1\frac{1}{x^{\frac{1}{p}-\frac{\epsilon}{p}}(1+x)}\,dx \Big)
\end{align*}
by integration by parts. From this we notice that
$$\epsilon I(\epsilon)\to\int_0^\infty \frac{1}{t^{\frac{1}{p}}(1+t)}\,dt=\frac{\pi}{\sin\frac{\pi}{p}}$$
when $\epsilon \to 0^+$. This together with \eqref{eq 13}, \eqref{eq 14}, \eqref{ineq 15} and \eqref{ineq 16} implies
\[
\frac{\sum\nolimits_{m,n=1}^\infty\big(\frac{n}{m}\big)^{\frac{1}{q}-\frac{1}{p}}\frac{a_m(\epsilon)b_n(\epsilon)}{m+n-1}}{\|(a_m(\epsilon))\|_{\ell^p}\|(b_n(\epsilon))\|_{\ell^q}}\geq\frac{\epsilon I(\epsilon)}{(1+\epsilon\,\phi(\epsilon))^{\frac{1}{p}} (1+\epsilon\,\psi(\epsilon))^{\frac{1}{q}}}\to\frac{\pi}{\sin\frac{\pi}{p}},
\]
when $\epsilon \to 0^+$.
\end{proof}

\section{The norm of the Hilbert matrix on the space $K^p$}

One can easily check that $\mathcal{H}$ induces a bounded operator on the space $K^p$, for $1<p<\infty$. Our second result is the determination of the exact value of the norm $\|\mathcal{H}\|_{K^p\to K^p}$. To that effect we shall use the variant of Hilbert's inequality in our Theorem \ref{theorem}. 

\begin{theorem}\label{Th1}
	If $1<p<\infty$, then
	\[
	\|\mathcal{H}\|_{K^p\to K^p}=\frac{\pi}{\sin\frac{\pi}{p}}
	\]
\end{theorem}
\begin{proof} 
	Let $f(z)=\sum\limits_{m=0}^\infty a_m z^m \in K^p$. Then
	\[
	\mathcal{H}(f)(z) = \sum_{n=0}^\infty \Big(\sum_{m=0}^\infty \frac{a_m}{m+n+1}\Big) z^n,
	\]
	and
	\begin{align*}
		\|\mathcal{H}(f)\|_{K^p}&=\Big(\sum_{n=0}^\infty(n+1)^{p-2} \Big|\sum_{m=0}^\infty \frac{a_m}{m+n+1}\Big|^p\Big)^{\frac 1p}\\
		&=\Big(\sum_{n=0}^\infty\Big|\sum_{m=0}^\infty (n+1)^{\frac{p-2}{p}}  \frac{a_m}{m+n+1}\Big|^p\Big)^{\frac 1p}.
	\end{align*}
Due to the duality of $\ell^p$ spaces 
	\[
	 \|\mathcal{H}(f)\|_{K^p}=\sup_{\|(b_n)\|_{\ell^q}=1}\Big|\sum_{m,n=0}^\infty(n+1)^{\frac{p-2}{p}}  \frac{a_mb_n}{m+n+1}\Big|,
	\]
where $\frac{1}{p}+\frac{1}{q}=1$.\newline 
Setting $A_m={a_m}{(m+1)^{\frac{p-2}{p}}}$, we have that $\|(A_m)\|_{\ell^p}=\|f\|_{K^p}$ and
\[
\sup_{\|f\|_{K^p}=1}\|\mathcal{H}(f)\|_{K^p}
=\sup_{\substack{ \|(A_m)\|_{\ell^p}=1,\\ \|(b_n)\|_{\ell^q}=1} }\, 
 \Big|\sum_{m,n=0}^\infty\Big(\frac{n+1}{m+1}\Big)^{\frac{1}{q}-\frac{1}{p}}\frac{A_mb_n}{m+n+1}\Big|\,=\,\frac{\pi}{\sin\frac{\pi}{p}}\,,
\]
because of Theorem \ref{theorem}.
\end{proof}


\bibliography{Literature}
\bibliographystyle{plain}

\end{document}